\documentclass[12pt]{article}

\input tcilatex

\input xy

\xyoption{poly}

\usepackage{amsfonts,epsf}
\usepackage{color}

\newcommand{\psdiag}[3]{\hspace{1mm}\raisebox{-#1mm}{\epsfysize#2mm

\epsffile{#3.eps}}\hspace{1mm}}

\begin{document}

\author{Rui Pedro Carpentier\\
\\ {\small rcarpent@math.ist.utl.pt} \\
{\small\it Departamento de Matem\'{a}tica} \\ {\small\it Centro de
An\'{a}lise Matem\'{a}tica, Geometria e Sistemas Din\^{a}micos}\\
{\small\it Instituto Superior T\'{e}cnico}\\ {\small\it Avenida
Rovisco Pais, 1049-001 Lisboa}\\ {\small\it Portugal}}
\title{Extending triangulations of the $2$-sphere to the $3$-disk preserving a $4$-coloring\footnote{Supported by Funda\c c\~ao para a Ci\^encia e a Tecnologia, project Quantum Topology, POCI/MAT/60352/2004 and PPCDT/MAT/60352/2004, and project New Geometry and Topology, PTDC/MAT/101503/2008.}}
\date{}
\maketitle

\newpage

\begin{abstract}
In this paper we prove that any triangulation of a $2$-dimensional sphere with a strict $4$-coloring on its vertices can seen as the boundary of a triangulation of a $3$-dimensional disk with the same vertices and preserving the $4$-coloring.
\end{abstract}

{\bf keywords:} Graph colorings, triangulations, $\Delta$-structures



\section{Introduction}

In \cite{8} it is stated that if a triangulation of a sphere can be colored on its vertices with exactly $4$ colors then this triangulation can be obtained from the tetrahedron by the following sequences of moves:

$$\psdiag{8}{16}{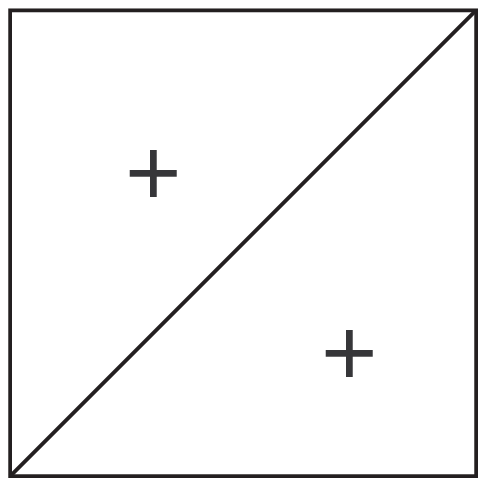} \stackrel{\mbox{move
I}}{\longleftrightarrow}\quad \psdiag{8}{16}{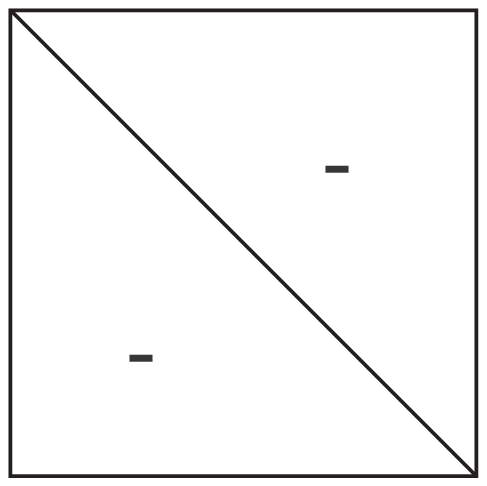}$$

$$\psdiag{8}{16}{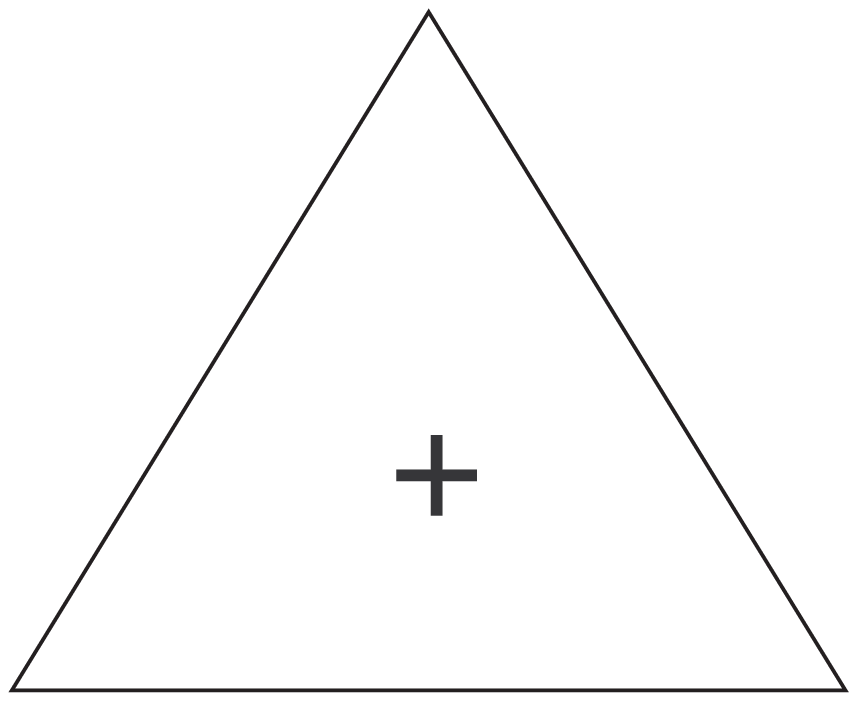} \stackrel{\mbox{move
II}}{\longrightarrow} \psdiag{8}{16}{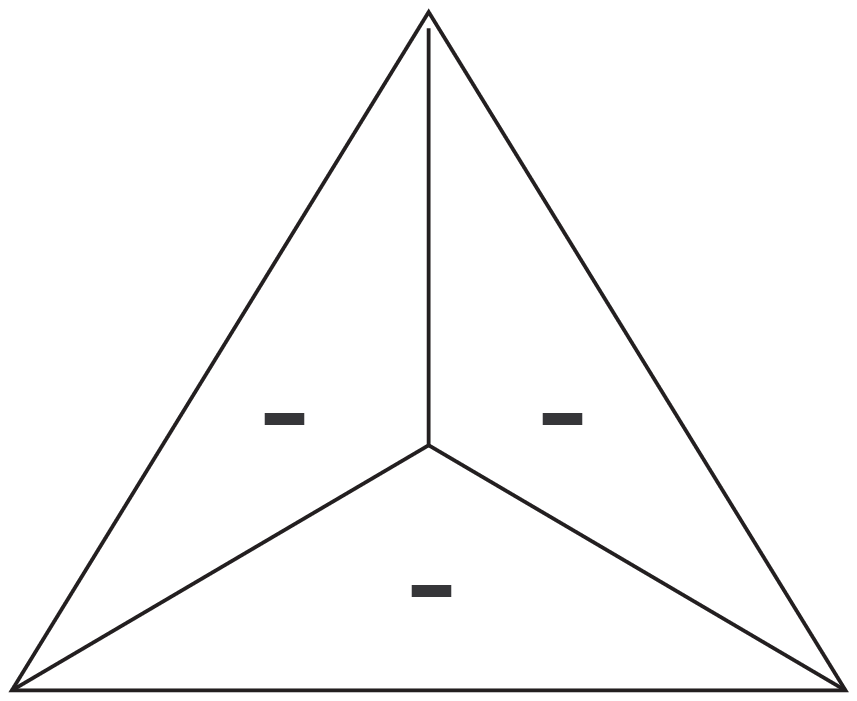}$$

$$\psdiag{8}{16}{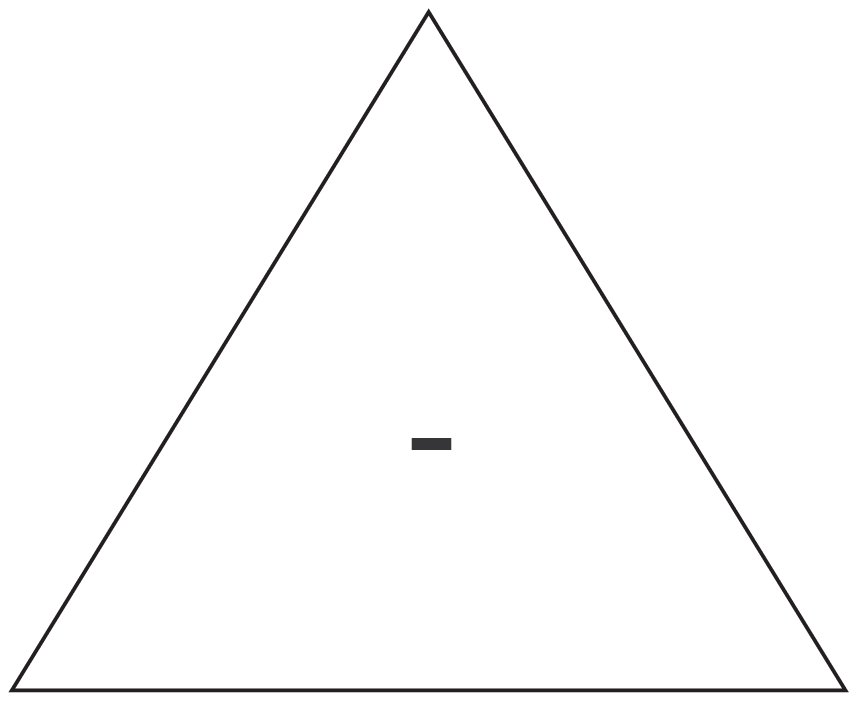} \stackrel{\mbox{move
II}}{\longrightarrow} \psdiag{8}{16}{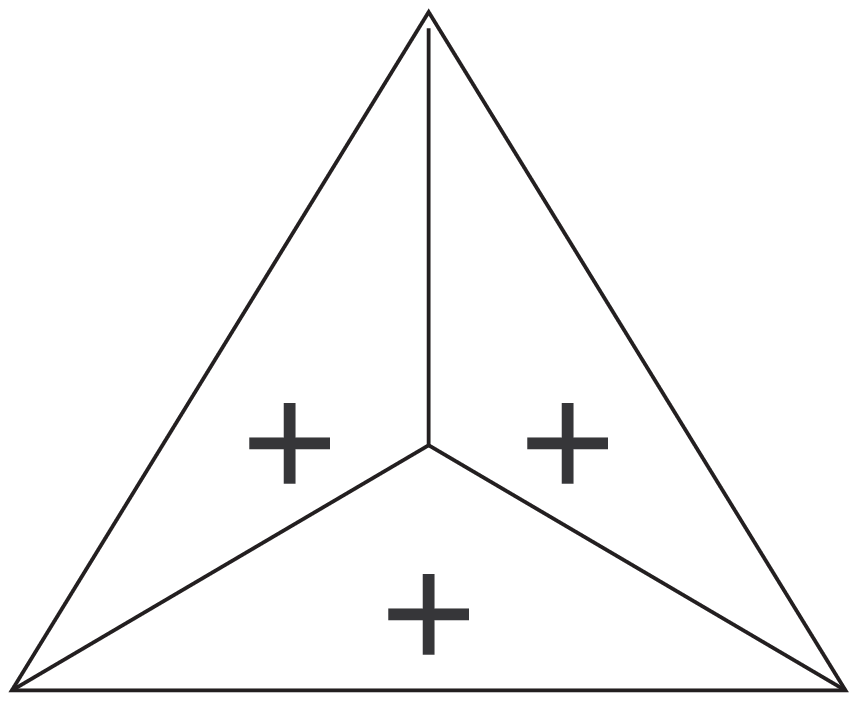}$$
where the signs are defined by the $4$-coloring in the following manner. A proper $4$-coloring on the
vertices of the triangulation induces a $3$-coloring on the edges of the
triangulation such that the boundary of each triangle is colored by
the three colors. This can be done by considering the four colors
as the four elements of the field $\mathbb{F}_4$ of order $4$, and then 
the $3$-coloring on the edges of the triangulation is obtained
by coloring each edge with the sum (or difference) of the colors
of its end points. By fixing an order on the three colors we get
a signing on the triangles, $+$ if the colors on the boundary are
ordered in the counterclockwise sense and $-$ if the colors on the
boundary are
ordered in the clockwise sense.
Note that move I does not change the coloring on the vertices and the moves II extend the coloring on the vertices in a unique way to the central vertex.

These moves are, disregarding the sign and unidirectional arrows, the Pachner moves \cite{2} which can be seen as the gluing of tetrahedra onto the triangulation.
This observation led the autor to wonder whether or not any triangulation of a strictly $4$-colorable sphere (i.e. one that can be colored with exactly $4$ colors), can be seen as the boundary of a triangulation of the $3$-dimensional disk with the same vertices and that preserves their colorings.

The answer is to this question in the affirmative is the main result of this paper, and the result of \cite{8} stated at the start of this introduction follows from this as a corollary.

\section{Triangulations, $\Delta$-structures and other algebraic topological prerequesites}

In this section we are going to clarify some topological definitions that we will use in this paper. In \cite{8}, Eliahou, Gravier and Payan defined a triangulation of a closed surface $S$ as a finite graph, loop-free but possibly with multiple edges, embedded in the surface $S$ and subdividing it into triangular faces. However, this definition, by allowing parallel edges, differs from the usual meaning given by many topologists which defines a triangulation of $S$ as an ordered pair $(K,h)$ where $K$ is a simplicial complex and $h:|K|\longrightarrow S$ is a homeomorphism from the geometric realization of $K$ to the space $S$ (e.g. see \cite{11}). Instead, we use as the definition of a triangulation on a space $X$ what is called in \cite{12} a $\Delta$-complex structure on $X$. 

Let $$\Delta^n=\{(t_0,\dots, t_n)\in\mathbb{R}^{n+1}:\sum_i t_i=1 \mbox{ and } t_i\geq 0 \mbox{ for all } i\}$$ be the {\it standard $n$-simplex}. A $k$-dimensional face of $\Delta^n$ (with $k<n$) is a subset of $\Delta^n$ where $n-k$ coordinates $t_{i_1},\dots, t_{i_{n-k}}$ are equal to zero. The union of all the faces of $\Delta^n$ is the {\it boundary} of $\Delta^n$, denoted by $\partial\Delta^n$. The {\it open $n$-simplex} $\mathring{\Delta}^n$ is $\Delta^n-\partial\Delta^n$, the interior of $\Delta^n$:
$$\mathring{\Delta}^n=\{(t_0,\dots, t_n)\in\mathbb{R}^{n+1}:\sum_i t_i=1 \mbox{ and } t_i>0 \mbox{ for all } i\}$$

A {\it $\Delta$-complex structure} on a space $X$ (which we shall also call a triangulation on $X$) is a (finite) collection of maps $\sigma_\alpha:\Delta^n\longrightarrow X$, with $n$ depending on the index $\alpha$, such that:

\begin{enumerate}
\item[(i)] The restriction $\sigma_\alpha|\mathring{\Delta}^n$ is injective, and each point $X$ is in the image of exactly one such restriction $\sigma_\alpha|\mathring{\Delta}^n$.

\item[(ii)] Each restriction $\sigma_\alpha$ to a face (of dimension $k<n$) of $\Delta^n$ is one of the maps $\sigma_\beta:\Delta^k\longrightarrow X$. Here we are identifying the ($k$-dimensional) face of $\Delta^n$ with $\Delta^k$ by the canonical linear homeomorphism between them that preserves the ordering of the vertices.

\item[(iii)] A set $A\subset X$ is open iff $\sigma_\alpha^{-1}(A)$ is open in $\Delta^n$ for each $\sigma_\alpha$.
\end{enumerate}

Given a space $X$ provided with a $\Delta$-complex structure $T$ we call it a {\it triangulated space}. 

We call the images of the $0$-simplex, $\sigma_\alpha:\Delta^0\longrightarrow X$ the {\it vertices} in $X$, the images of the $1$-simplex, $\sigma_\alpha:\Delta^1\longrightarrow X$ the {\it edges} in $X$, the images of the $2$-simplex the {\it triangles} in $X$, and the images of the $3$-simplex the {\it tetrahedra} in $X$. The {\it graph} of a triangulation is formed by its sets of vertices and edges, and a (proper) coloring on a triangulation is a (proper) coloring on its graph.

The {\it star} of a vertex $v$ in $X$ is formed by the images of the maps $\sigma_\alpha:\Delta^n\longrightarrow X$ that contain $v$, $$\mbox{st}(v)=\bigcup_{v\in\sigma_\alpha(\Delta^n)}\sigma_\alpha(\Delta^n)$$

The {\it deletion} of a vertex $v$ in $X$ is formed by the images of the maps $\sigma_\alpha:\Delta^n\longrightarrow X$ that do not contain $v$, $$\mbox{dl}(v)=\bigcup_{v\not\in\sigma_\alpha(\Delta^n)}\sigma_\alpha(\Delta^n)$$

A $\Delta$-complex structure on a space $X$ induces a $\Delta$-complex structure on its cone $CX=X\times [0,1]/X\times \{0\}$ in a natural way: for each map $\sigma_\alpha :\Delta^k \longrightarrow X$ we have two maps $\overline{\sigma}_\alpha :\Delta^k \longrightarrow CX$ given by $\overline{\sigma}_\alpha(t_0,\dots,t_k)=(\sigma_\alpha(t_0,\dots,t_k),1)$ and $\hat{\sigma}_\alpha :\Delta^{k+1} \longrightarrow CX$ given by $$\hat{\sigma}_\alpha(t_0,\dots,t_{k+1})=\left(\sigma_\alpha(\frac{t_0}{1-t_{k+1}},\dots,\frac{t_k}{1-t_{k+1}}),1-t_{k+1}\right)$$ (for $t_{k+1}\not=1$) and $\hat{\sigma}_\alpha(0,\dots,0,1)=*$ where $*$ is the point $X\times\{0\}$ in $CX$. Finally we complete the triangulation with the map $\sigma_* :\Delta^0 \longrightarrow CX$ given by $\sigma_*(\Delta^0)=*$. We call this construction the {\it cone of the triangulation} of $X$.


\section{The main result and a corollary}

Our main result can be stated as follows:

\begin{theorem}

If $\psi$ is a strict $4$-coloring of a triangulation $T$ of
the $2$-sphere $\mathbb{S}^2$ then there exists a triangulation $T'$
of the $3$-disk $\mathbb{D}^3$ such that {$T$ is the triangulation induced by $T'$ on the boundary of the disk}, the vertices
of $T'$ are in $T$ and $\psi$ is a $4$-coloring of $T'$.
\label{5}
\end{theorem}

To prove this theorem we will make use of the following lemma, which can be seen as a version of the previous theorem for one dimension lower:

\begin{lemma}

If $\psi$ is a strict $3$-coloring of a triangulation $T$ of
the circle $\mathbb{S}^1$ then there exists a triangulation $T'$
of the disk $\mathbb{D}^2$ such that {$T$ is the triangulation induced by $T'$ on the boundary of the disk}, the vertices
of $T'$ are in $T$ and $\psi$ is a $3$-coloring of $T'$.
\label{6}
\end{lemma}

\TeXButton{Proof}{\proof}

In this case $T$ is just a cycle graph, thus if $\psi$ is a
strict $3$-coloring then there are three consecutive vertices
$v_1$, $v_2$ and $v_3$ with distinct colors $a$, $b$ and $c$.
Suppose that the middle vertex $v_2$ is colored by $b$. If $v_2$ is
the only vertex in $T$ colored by $b$ then we can add edges
linking $v_2$ with all vertices of $T$ and therefore we get a
triangulation $T'$ of the disk $\mathbb{D}^2$ with the desired
properties.

$$\psdiag{10}{20}{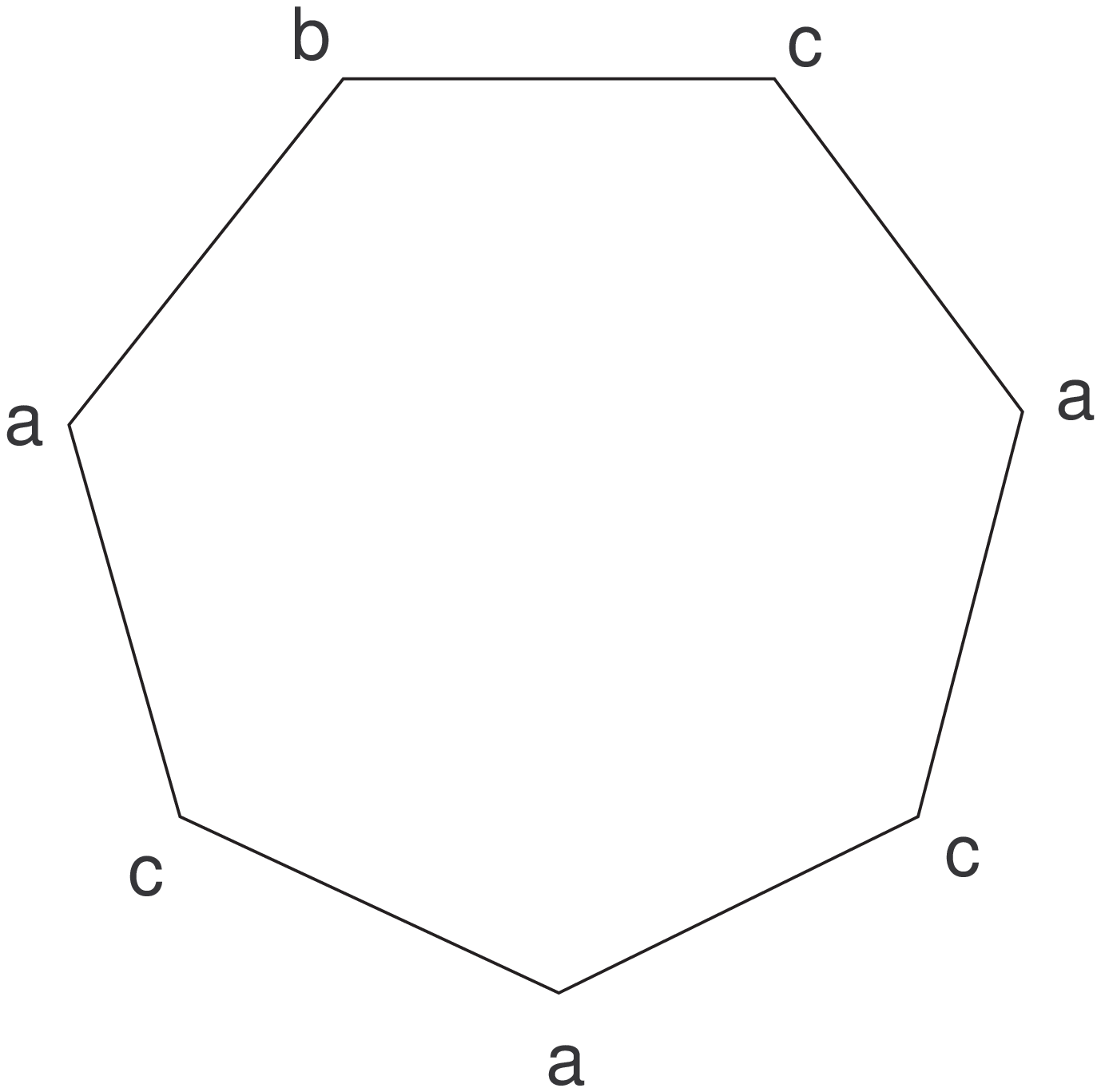} \longrightarrow \psdiag{10}{20}{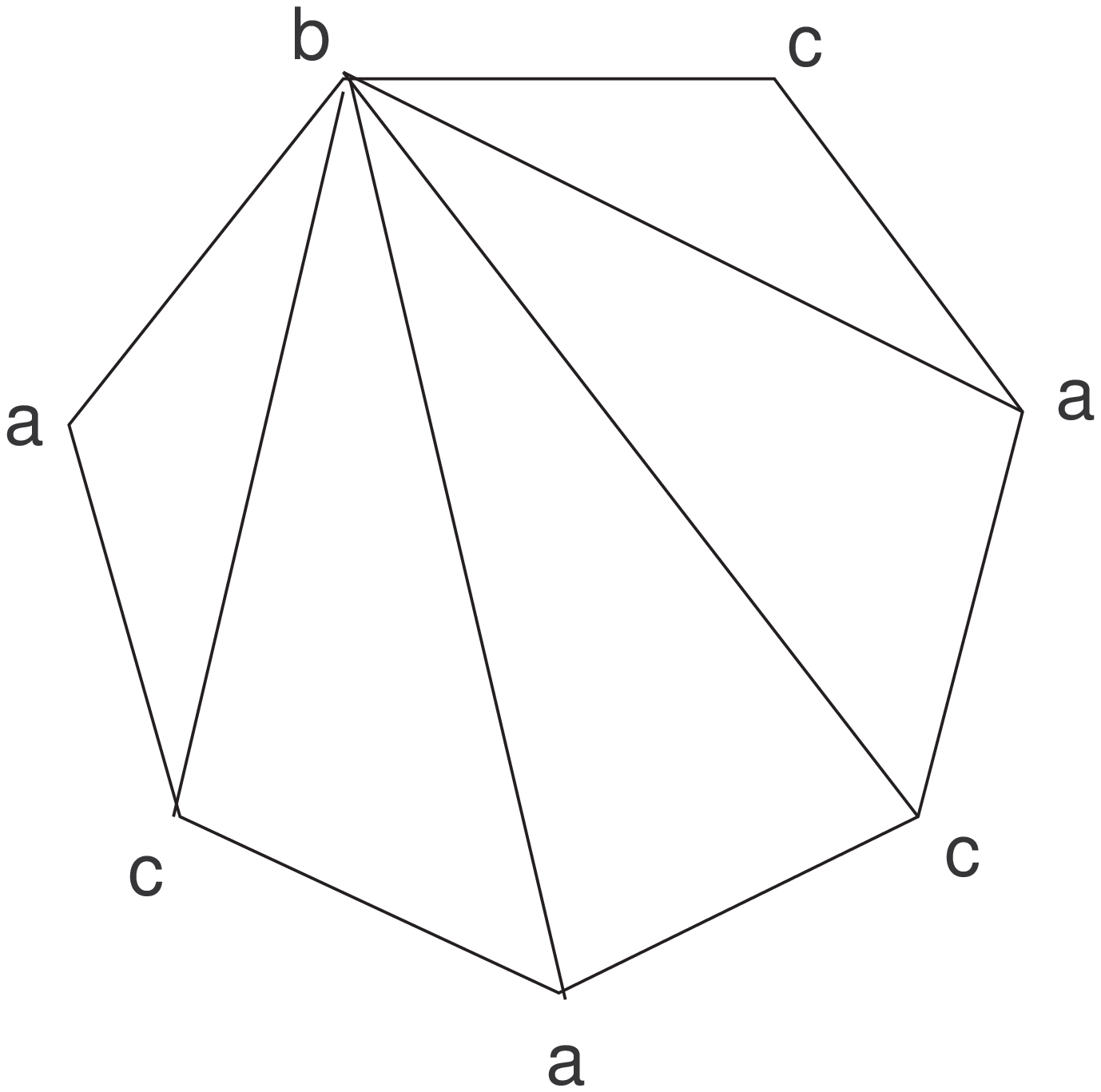}$$

If $v_2$ is not the only vertex in $T$ colored by $b$ then we
can add an edge linking $v_1$ with $v_3$, and then complete, by
induction {on the number of vertices of $T$}, the triangulation on the disk whose boundary is the cycle $v_1$, $v_3$, ...
,$v_n$.

$$\psdiag{10}{20}{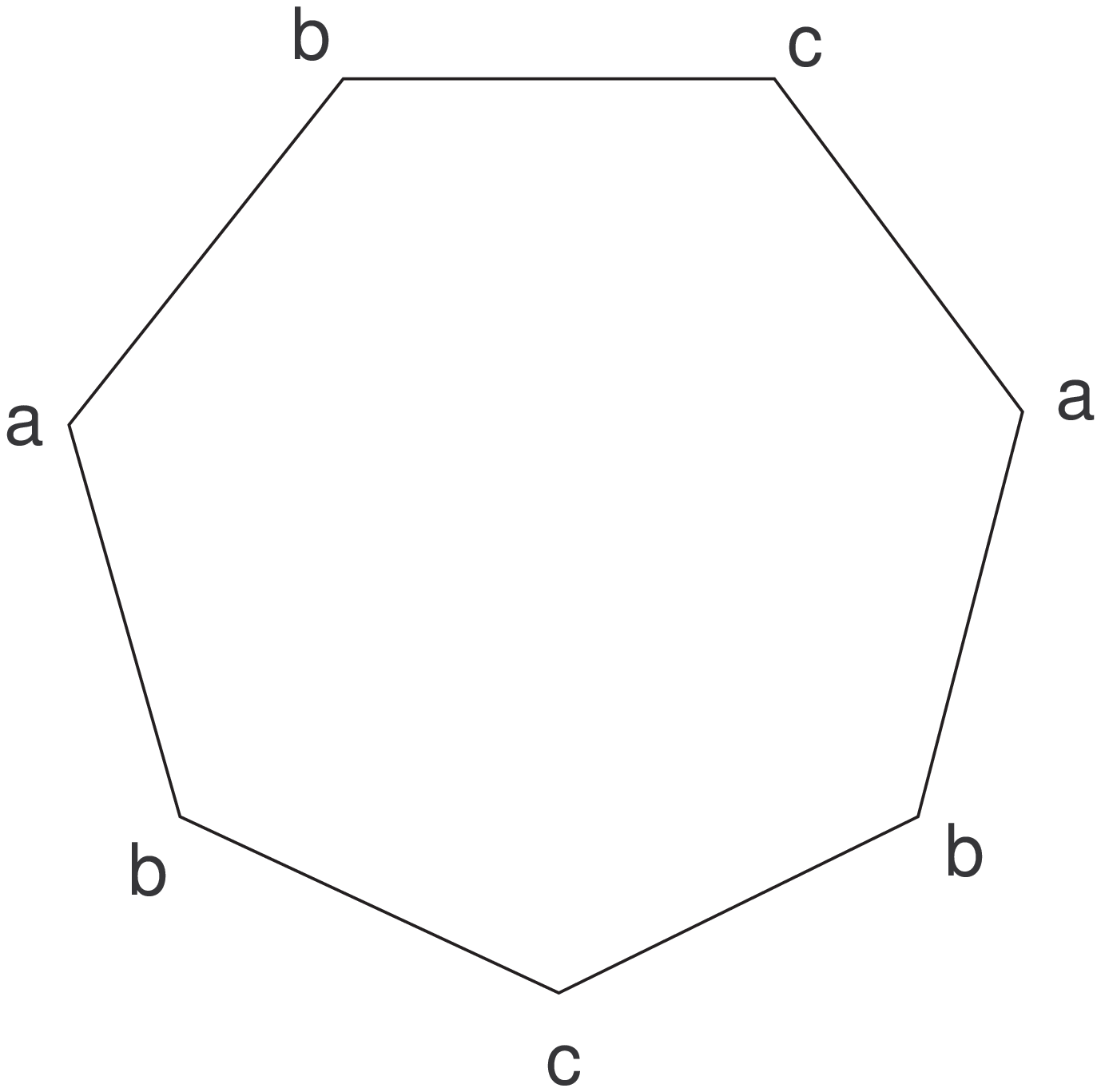} \longrightarrow
\psdiag{10}{20}{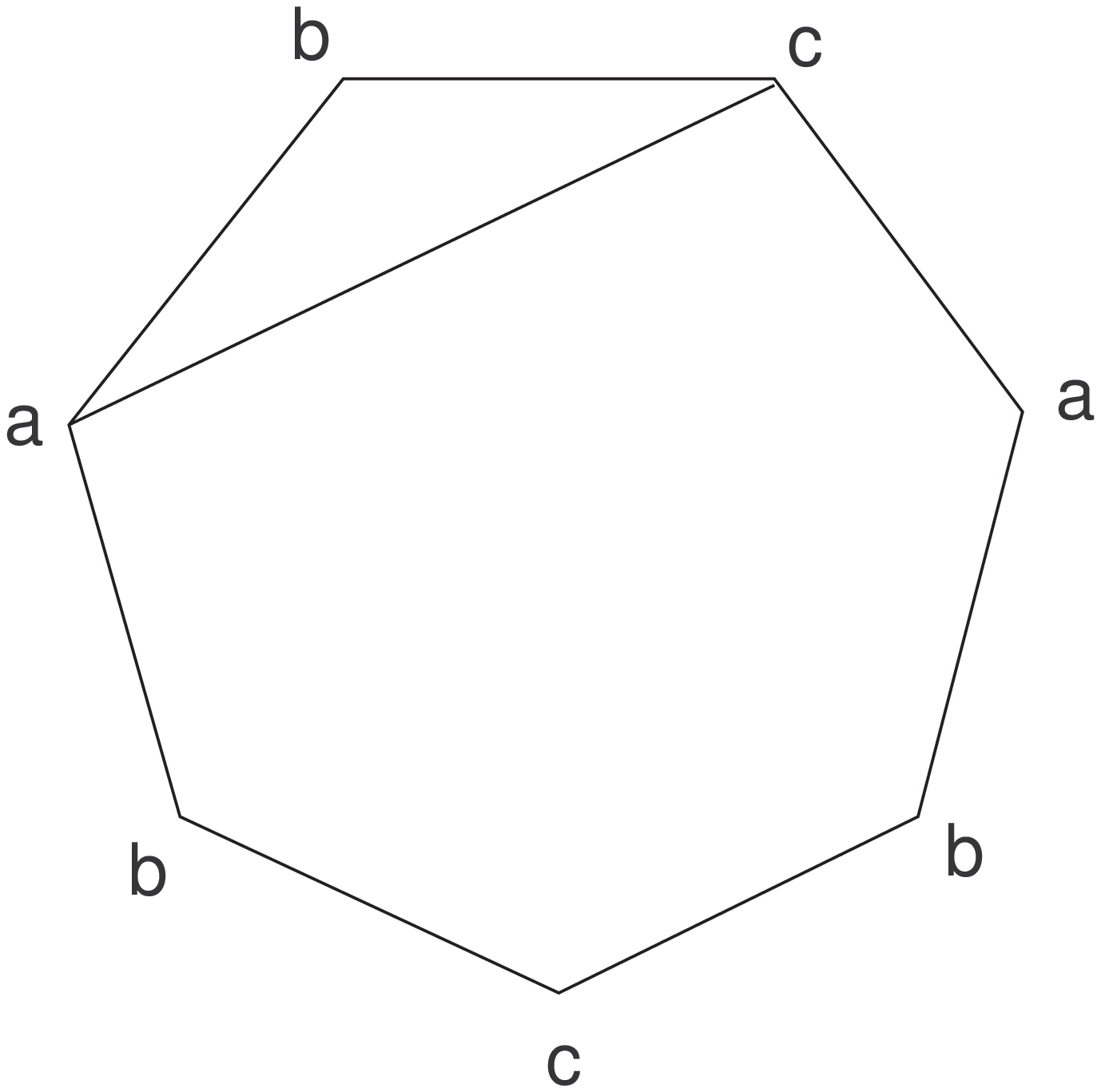}\stackrel{\mbox{induction}}{\longrightarrow}
\psdiag{10}{20}{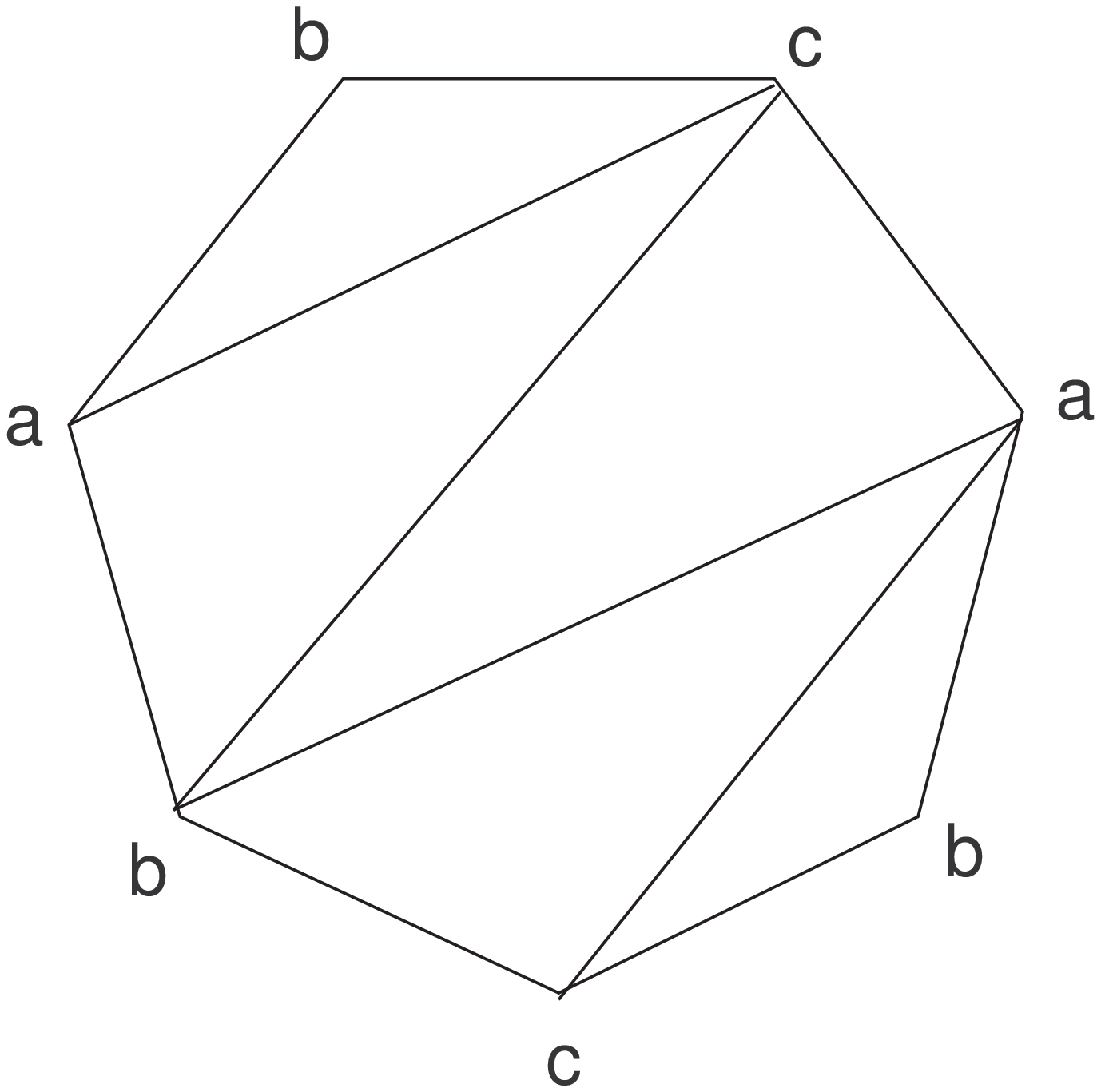}$$

\TeXButton{End Proof}{\endproof}

{\bf Proof of theorem \ref{5}.}

First we consider the case where the triangulation $T$ has no parallel edges (i.e. multiple edges between the same pair of vertices), this means that the star and the deletion (see section 2) of any vertex of $T$ are simplicial disks.

We start by proving that in the triangulation $T$ there exists a vertex $v$ that is adjacent to a cycle colored by three colors. We take a triangle colored by (say)
$a$, $b$ and $c$ and consider the region formed by the triangles colored by the same colors.
Since $\psi$ is a strict $4$-coloring, this region has a non-empty boundary and any vertex of
the boundary satisfies the required condition, because its link necessarily has vertices colored by $d$ and by two colors amongst $a$, $b$ and $c$ which differ from the color of the vertex itself.

Now, if $v$ is the only vertex in $T$ colored by its color $\psi(v)$ then we take the cone of the deletion of $v$ to
 get a triangulation $T'$ of the disk
$\mathbb{D}^3$ with the desired properties.

If, on the contrary, $v$ is not the only vertex in $T$ colored by $\psi(v)$, we remove $v$, use the lemma to triangulate the
region $R$ bounded by the link of $v$ to get a triangulation of $\mathbb{S}^2$, $T''$, with one less vertex than $T$ which is still strict $4$-colored, use induction to get a
triangulation of the disk $\mathbb{D}^3$ bounded by this and finally attach to the region $R$ its cone to
obtain the desired triangulation $T'$ of $\mathbb{D}^3$.

In the case when the triangulation does have parallel edges, we take two parallel edges $e$ and $e'$ linking two vertices $v_i$ and $v_j$ colored by two colors (say $1$ and $2$). The two parallel edges form a cycle $C$ that splits the sphere into two (triangulated) $2$-disks $D_1$ and $D_2$.

First, we suppose that we have a strict $4$-coloring for at least one of the disks (say $D_1$ has vertices colored by all four colors). Thus glueing the edges $e$ and $e'$ we get a (triangulated) sphere $S_1$ with less vertices\footnote{Note that since the cycle $C$ is formed by only two edges each disk must contain inner vertices.} than the original sphere. To the other disk we attach along the cycle $C$ its cone in order to get another (triangulated) sphere $S_2$ with less vertices\footnote{Note that since the disk $D_1$ has vertices of all colors at least two vertices must be deleted from the original triangulated sphere in order to get the disk $D_2$.} than the original sphere. {The disk $D_2$ is colored at least by three colors (say $1$, $2$ and $3$), so we color the inserted vertex on $S_2$ with the color $4$ in order to guarantee we have a strict $4$-coloring on $S_2$}. Therefore, by induction on the number of vertices, $S_1$ is the boundary of a triangulated $3$-disk ${\cal D}'_1$ without inner vertices and preserving the $4$-coloring and $S_2$ is the boundary of a triangulated $3$-disk ${\cal D}'_2$ with the same properties.

Now we perform the following surgery: take the edge $e^*$ on $S_1$ resulting from the glueing of the parallel edges $e$ and $e'$ and search in ${\cal D}'_1$ for a triangle $t$ containing that edge and a vertex colored with the color $4$ (such a triangle exists because the edge belongs to some tetrahedron). 

If $t$ is an inner triangle in ${\cal D}'_1$ adjacent to two tetrahedra $\tau_1$ and $\tau_2$, we replace $t$ by the copies of it $t_1$ (adjacent to $\tau_1$) and $t_2$ (adjacent to $\tau_2$) sharing the same vertices and the same edges except the edge $e^*$ which is replaced by the old parallel edges $e$ and $e'$ (in other words, we are opening a fissure in ${\cal D}'_1$ through the edge $e^*$ along the triangle $t$), then we identify the triangles $t_1$ and $t_2$ with the triangles in $S_2$ (which are in ${\cal D}'_2$) that made the cone of the cycle $C$ (formed by $e$ and $e'$). This ``grafting'' of the disk ${\cal D}'_2$ in the disk ${\cal D}'_1$ produces the desired triangulated $3$-disk.


{If $t$ is on the surface $S_1$ then there exists a triangle $t'$ on $D_1$ with the same vertices. If the triangle $t'$ is adjacent to the edge $e$ (resp. $e'$) then we identify the triangle $t$ with the triangle on $S_2$ that belongs to the cone of the cycle $C$ and is adjacent to the edge $e'$ (resp. $e$).
This glueing of the disk ${\cal D}'_1$ with the disk ${\cal D}'_2$ produces the desired triangulated $3$-disk.}


Finally, if neither of the disks $D_1$ and $D_2$ has a strict $4$-coloring (we can suppose that $D_1$ is colored by $1$, $2$ and $3$ and $D_2$ is colored by $1$, $2$ and $4$) then we remove from $D_1$ the edge $e$ and the triangle $t_1$ in $D_1$ incident to it and we remove from $D_2$ the edge $e'$ and the triangle $t_2$ in $D_2$ incident to it. We then get two new $2$-disks $D'_1$ and $D'_2$. We take the triangulated cones of them both. Let $v_1$ be the opposite vertex to the edge $e$ in the triangle $t_1$, $v_2$ the opposite vertex to the edge $e'$ in the triangle $t_2$, $v^*_1$ the cone vertex of the cone of $D'_1$ and $v^*_2$ the cone vertex of the cone of $D'_2$. Then we get the desired triangulation of the $3$-disk by glueing the two cones by identifying the triangle in the cone of $D'_1$ with vertices $v_1$, $v_i$ and $v^*_1$ with the triangle in the cone of $D'_2$ with vertices $v_2$, $v_i$ and $v^*_2$ and the triangle in the cone of $D'_1$ with vertices $v_1$, $v_j$ and $v^*_1$ with the triangle in the cone of $D'_2$ with vertices $v_2$, $v_j$ and $v^*_2$ (recall that the two ends of the parallel edges $e$ and $e'$ were called $v_i$ and $v_j$, and note that, through the identifications, the triangles $t_1$ and $t_2$ appear on the boundary of the 3-disk as faces of the cones of $D_2'$ and $D_1'$ respectively.).
\TeXButton{End Proof}{\endproof}

As a corollary we have an alternative Proof for Theorem 1.3 of \cite{8} which may be stated as follows:

\begin{theorem}\label{4}

Suppose we are given a triangulation $T$ of the sphere with signed faces.
Then the signing comes from a strict $4$-coloring of $T$ (i.e. a
$4$-coloring that uses the four colors) if and only if $T$ comes
from the tetrahedron with the same sign on all its faces, by means of a
sequence of signed diagonal flips (move I) and$/$or divisions of a
triangle into three triangles (by adding a vertex $v$ inside the
triangle and edges joining $v$ to the vertices of the triangle)
with opposite signs (move II).

\end{theorem}

$$\psdiag{8}{16}{creassig1p} \stackrel{\mbox{move
I}}{\longleftrightarrow}\quad \psdiag{8}{16}{creassig2p}$$

$$\psdiag{8}{16}{moveII1} \stackrel{\mbox{move
II}}{\longrightarrow} \psdiag{8}{16}{moveII2}$$

$$\psdiag{8}{16}{moveII3} \stackrel{\mbox{move
II}}{\longrightarrow} \psdiag{8}{16}{moveII4}$$

\TeXButton{Proof}{\proof}

For the sufficient condition the proof is the same as the proof in
\cite{8}. Move I does not change the coloring and move II
extends the coloring in a unique way.

For the other implication, we take the triangulation of $\mathbb{D}^3$ obtained in Theorem \ref{5},
we choose one tetrahedron and by adding the adjacent tetrahedra one by one we get a sequence
of moves I and II (the signs of the faces are determined by the coloring as was observed in
the paragraph following theorem 1).

We only have to see that this sequence of attaching adjacent tetrahedra can be done keeping the topology of a $3$-disk in each step. We are going to prove by induction on the number of tetrahedra that this can be done independently of the tetrahedron we choose to start.

First we note that the triangulation produced in Theorem \ref{5} is obtained by one of the following procedures: 1) taking the triangulated cone of a triangulated $2$-disk, 2) attaching a triangulated cone of a triangulated $2$-disk without inner vertices to a previously triangulated $3$-disk, 3) ``grafting'' a triangulated $3$-disk into another triangulated $3$-disk, 4) glueing two previously triangulated $3$-disks along a shared triangle on their surface,  and 5) glueing two triangulated cones of $2$-disks along two adjacent shared triangles on their surface.

In the first procedure we have to prove that, given a triangulation of the $2$-disk and a triangle in it, there exists a sequence of attaching adjacent triangles, starting from the given triangle and ending with the given triangulation, such that in each step we have a triangulation of a $2$-disk. This can be easily proved by induction on the number of triangles of the triangulation. If there exists a edge that cuts the $2$-disk in two then we get, by induction, a sequence of attaching triangles on the disk that contains the starting triangle and we continue the sequence on the other disk starting with the triangle incident to the cutting edge. If there are no such cutting edges and the triangulation has more than one triangle (the case of a one-triangle triangulation is trivial) then we choose one triangle incident to a boundary edge different from the starting triangle, we remove it and use induction in order to get a sequence of attaching triangles from the starting triangle to this triangulation of the disk with the triangle removed and then we complete the sequence by attaching this last triangle.

In the second procedure, we have by induction sequences for any starting tetrahedron both in the smaller $3$-disk and in the attaching cone. If the starting tetrahedron is in the smaller $3$-disk then we take a sequence in the smaller $3$-disk for that tetrahedron and complete it with a sequence in the attaching cone. If the starting tetrahedon is in the cone then take a sequence in the cone starting at that tetrahedron and a sequence in the smaller $3$-disk starting at the tetrahedron adjacent to that tetrahedron, and we proceed in the following way. Start with the starting tetrahedron, follow it with the sequence in the smaller $3$-disk, and complete the sequence with the rest of the sequence in the cone.

In the third ``grafting'' case, if the starting tetrahedron is in the $3$-disk ${\cal D}_2'$  that is grafted into the other $3$-disk ${\cal D}_1'$ (see the final part of the proof of theorem \ref{5}) using an inner triangle, then we take, by induction, a sequence in ${\cal D}_1'$  starting at that tetrahedron and complete the sequence in the other disk by starting at one of the tetrahedra incident to the triangle where the ``grafting'' is done. If the starting tetrahedron is in the other $3$-disk (${\cal D}_1'$) then from a sequence in that disk starting at that tetrahedron, given by induction, we take the partial sequence from that starting tetrahedron to the first tetrahedron in the sequence that is incident to the triangle where the ``grafting'' is done, then we follow with a sequence in the grafted disk ${\cal D}_2'$ that starts at the tetrahedron adjacent to the latter tetrahedron, and finally we complete with the rest of the sequence in ${\cal D}_1'$ .

In the fourth case of two $3$-disks attached by a single triangle, we start with a given tetrahedron, take the sequence given by induction in the disk that contains it, and complete the sequence in the other disk starting at the tetrahedron incident to the glueing triangle.

In the final case, we have the cones of two disks $D'_1$ and $D'_2$ glued by identifying a certain pair of adjacent triangles on each cone. Without loss of generality suppose that the starting tetrahedron is in the cone of the disk $D'_1$. We know that for a cone of a triangulated $2$-disk there is a sequence of attaching tetrahedra from any starting tetrahedron to the whole triangulation keeping the $3$-disk topology in each step. So we can find such a sequence for the cone of the disk $D'_1$ starting at the given tetrahedron, say $T_1, T_2, \cdots, T_k$, and another such sequence for the cone of the disk $D_2$ (i.e. the disk $D'_2$ with the triangle with vertices $v_i$, $v_2$ and $v_j$ attached to it - see the last part of the proof of Theorem \ref{5}), starting at the cone of this last triangle, say $T_1', T_2', \cdots , T_l'$. Thus we obtain a sequence for the whole 3-disk by omitting $T_1'$ and composing these two sequences as $T_1, T_2, \cdots, T_k, T_2', \cdots , T_l'$.
\TeXButton{End Proof}{\endproof}

\section{Comments and a conjecture}

We have seen the analogy between Theorem~\ref{5} and Lemma~\ref{6} and how the second result is used in the proof of the first. This leads us to conjecture the following:

\begin{conjecture}\label{c}
If $\psi$ is a strict $n+2$-coloring of a triangulation $T$ of
the $n$-sphere $\mathbb{S}^n$ then there exists a triangulation
$T'$ of the disk $\mathbb{D}^{n+1}$ such that $\partial T'=T$, the
vertices of $T'$ are in $T$ and $\psi$ is a $n+2$-coloring of
$T'$. 
\end{conjecture}

It is not clear if and how the proof of the Theorem~\ref{5} can be adapted to higher dimensions even in the weaker case of triangulations without parallel edges (the usual definition for triangulation). This because, for triagulations of spheres of dimension greater than $4$, the link of a vertex is not necessarily a sphere of lower dimension. 
However, if we consider only the case of piecewise-linear triangulations (where the link of any simplex is a piecewise-linear sphere, see \cite{1}) the first part of the proof of Theorem~\ref{5} seems to apply recursively to prove the previous conjecture.
We also mention that, for piecewise-linear triangulations, a weaker version of this conjecture (by allowing inner vertices) follows from lemma 3.1 of \cite{3}.

Another open problem is to see to what extent the proof of Theorem~\ref{4} depends on the triangulations of the $3$-disk obtained in the proof of Theorem~\ref{5}. In other words, given a triangulation of the $3$-disk, not necessarily  of the type used in the proof of Theorem~\ref{5}, we want to know if there is a sequence of attaching tetrahedra from an initial tetrahedron to the final triangulation such that the topology of the $3$-disk is kept in each step.


{\bf Acknowledgment} - I wish to thank Roger Picken for his useful suggestions and comments. I also want to thank Nathalie Wahl for telling me about the relation between her result in \cite{3} and Conjecture~\ref{c}.

\end{document}